\documentclass[12pt]{amsart}
\usepackage{amsmath}
\usepackage{amsfonts,amssymb,amsthm}

\newcommand{\vol}{{\mathrm{vol}\,}}

\begin{document}


\title[An Extended Correction]{An extended correction to ``Combinatorial Scalar Curvature and Rigidity of Ball Packings,'' (by
D.~Cooper and I.~Rivin).}

\author{Igor Rivin}


\address{Department of Mathematics, Temple University, Philadelphia}

\curraddr{Mathematics Department, Princeton University}

\email{rivin@math.temple.edu}

\thanks{The author would like to thank David Glickenstein and Bennett Chow (both of UCSD Mathematics Department)
for bringing the error to his attention.}

\date{today}



\begin{abstract}
It has been pointed out to the author by David Glickenstein that
the proof of the (closely related) Lemmas 1.2 and 3.2 in
\cite{balls} is incorrect. The statements of both Lemmas are
correct, and the purpose of this note is to give a correct
argument. The argument is of some interest in its own right.
\end{abstract}

\maketitle

\section*{Introduction}
Let us first recall the setup of \cite{balls}. In that paper we
study \emph{conformal} simplices. These are simplices $T(r_1, r_2,
r_3, r_4)$ in $3$-dimensional spaces of constant curvature such
that there are positive numbers $r_1, r_2, r_3, r_4,$ such that
the length $l_{ij}$ of the edge joining the $i-$th and the $j-$th
vertex of the simplex is given by $l_{ij} = r_i+r_j.$

On the set of conformal simplices we define a function $S,$ as
follows:
\begin{equation*}
S(r_1, r_2, r_3, r_4) = \begin{cases} \sum_{i=1}^4 r_i S_i \qquad
\mbox{for simplices in $\mathbb{E}^3,$}\\
2 \vol T + \sum_{i=1}^4 r_i S_i \qquad \mbox{for simplices in
$\mathbb{H}^3,$}
\end{cases}
\end{equation*}
where $\vol$ stands for the hyperbolic volume of the simplex, and
$S_i$ stands for the solid angle at the $i$-th vertex: if
$\alpha_{ij}$ is the dihedral angle at the edge joining the $i-$th
and the $j-$th vertices, then
\begin{equation}
\label{solid}
 S_i = - \pi + \sum_{j\neq i} \alpha_{ij}.
\end{equation}
The key property of the function $S$ is, as shown in \cite{balls},
that
\begin{equation}
\label{hess} H(S)_{ij} =  \frac{\partial^2 S}{\partial r_i
\partial r_j} = \frac{\partial S_i}{\partial r_j},
\end{equation}
where we use $H(S)$ to denote the Hessian matrix of $S.$
Lemma 1.2 states that $H(S)$ is negative semi-definite for
simplices in $\mathbb{E}^3,$ with the zero direction spanned by
the vector $(r_1, r_2, r_3, r_4),$ and corresponding to the
rescaling deformation of the simplex. Lemma 3.2 states that $H(S)$
is negative definite for simplices in $\mathbb{H}^3.$
\section{Proofs of the Lemmas}
The proofs given in \cite{balls} work without modification when
all the radii are equal ($r_1 = r_2 = r_3 = r_4.$) Since the set
of all conformal simplices is connected (as shown in
\cite{balls}), it suffices to show that the rank of $H(S)$ always
equals $3$ in the Euclidean case and $4$ in the hyperbolic case.
The proof will rest on the following observations:

\begin{itemize}

\item{(a)} Define a function $R$ on the set of all simplices by
\[
R(l_{12}, l_{13}, l_{14}, l_{23}, l_{24},  l_{34}) = 
\begin{cases}
\sum_{i< j} l_{ij} (\pi - \alpha_{ij}),\qquad \mbox{Euclidean,}\\
\sum_{i< j} l_{ij} (\pi - \alpha_{ij}) + 2 \vol T ,\qquad \mbox{hyperbolic.}
\end{cases}
\]
and define the map 
\[
i(r_1, r_2, r_3, r_4) = (r_1 + r_2, r_1 + r_3, r_1 + r_4, r_2 +
r_3, r_2 + r_4, r_3 + r_4).
\]
As noted in \cite{balls}, 
\begin{equation}
\label{pullback}
S = i^*R.
\end{equation}
 Also as noted in \cite{balls},
\begin{equation}
\label{schlaf}
\frac{\partial R}{\partial l_{ij}} = \alpha_{ij},
\end{equation}
and so
\begin{equation}
\label{schlaf2}
\frac{\partial^2 R}{\partial l_{ij}\partial l_{kl}} = \frac{\partial
\alpha_{ij}}{\partial l_{kl}}.
\end{equation}
\item{(b)} Since simplices are infinitesimally rigid, the Hessian
matrix of $R$ is nonsingular in the hyperbolic case, and nonsingular
when restricted to the cone $\sum_{i<j} l_{ij} = 1,\quad l_{ij} > 0$
in the Euclidean case (the implication is given by
Eq.~(\ref{schlaf2}).) In particular, the Hessian matrix is nonsingular
when restricted to the image of $i.$ By Eq.~(\ref{pullback}) and
linearity of the map $i,$ it follows
that $S$ is nonsingular (on the cone $\sum r_i = 1$, in the Euclidean
case).
\end{itemize} It follows that to check the convexity of $S,$ we need only
do it at one point. As pointed out above, the argument given in
\cite{balls} works without modification for regular simplices, and so
the convexity of $S$ follows.

\bibliographystyle{amsplain}

\end{document}